\documentclass[12pt,a4paper,reqno]{amsart}
\usepackage{amssymb}
\usepackage{latexsym}
\usepackage{exscale}

\numberwithin{equation}{section}
\newtheorem{theor}{Theorem}[section]
\newtheorem{lemma}[theor]{Lemma}

\newtheorem{corol}[theor]{Corollary}

\newcommand{\RR}{\mathbb{R}}

\newcommand{\W}{\mathcal{W}}

\newcommand{\T}{\mathcal{T}}

\newcommand{\A}{\mathcal{A}}

\newcommand{\BMO}{{\rm BMO}}

\newcommand{\aver}[1]{-\hskip-0.46cm\int_{#1}}

\begin{document}
\allowdisplaybreaks

\title[Riesz transforms on manifolds and weights]{Weighted norm inequalities,  off-diagonal estimates and
elliptic operators.
\\[.2cm]
{\footnotesize Part IV: Riesz transforms on manifolds and weights}}

\author{Pascal Auscher}

\address{Pascal Auscher
\\
Universit\'e de Paris-Sud et CNRS UMR 8628
\\
91405 Orsay Cedex, France} \email{pascal.auscher@math.u-psud.fr}

\author{Jos\'e Mar{\'\i}a Martell}

\address{Jos\'e Mar{\'\i}a Martell
\\
Instituto de Matem\'aticas y F{\'\i}sica Fundamental
\\
Consejo Superior de Investigaciones Cient{\'\i}ficas
\\
C/ Serrano 123
\\
28006 Madrid, Spain}

\address{\null\vskip-.7cm and\vskip-.7cm\null}

\address{Departamento de Matem\'aticas \\ Universidad Aut\'onoma de Madrid \\
28049 Madrid, Spain } \email{chema.martell@uam.es}

\thanks{This work was partially supported by the European Union
(IHP Network ``Harmonic Analysis and Related Problems'' 2002-2006,
Contract HPRN-CT-2001-00273-HARP). The second author was also
supported by MEC ``Programa Ram\'on y Cajal, 2005'' and by MEC Grant
MTM2007-60952.}
\thanks{We warmly thank T. Coulhon for interest and help in finding appropriate references. We also want to express our deep gratitude to F. Bernicot and J. Zhao for letting us use their unpublished work. This have lead us to remove all use of Poincar\'e inequalities in the revised version. We thank the referee for suggestions to enhance the presentation of this article.
}

\date{December 29, 2006. \textit{Revised}: November 8, 2007}



\begin{abstract}
This is the fourth article of our series. Here, we  study weighted norm inequalities for the Riesz
transform  of the Laplace-Beltrami operator on Riemannian manifolds
and of subelliptic sum of squares on Lie groups, under   the doubling
volume property and Gaussian upper bounds.
\end{abstract}

\maketitle

\section{Introduction and main results}
On $\RR^n$,  it is well-known that the classical Riesz transforms
$R_{j}$, $1\le j \le n$, are bounded on $L^p(\RR^n,dx)$ for
$1<p<\infty$ and are of weak-type (1,1) with respect to $dx$. As a
consequence of the weighted theory for classical Calder\'on-Zygmund
operators, the Riesz transforms are also bounded on $L^p(\RR^n,
w(x)dx)$ for all $w\in A_p(dx)$, $1<p<\infty$, and are of weak-type
(1,1) with respect to $w(x)dx $ for $w\in A_{1}(dx)$. Furthermore, it
can be shown that the $A_{p}$ condition on the weight is necessary
for the weighted $L^p$ boundedness of the Riesz transforms (see, for
example, \cite{Gra}).

On a manifold, there has been a number of works discussing the
validity of the unweighted $L^p$ theory depending on the geometry of
the manifold.  Although some progress has been done in this
direction, the general picture is far from clear. A  difficulty is
that one has to leave the class of Calder\'on-Zygmund operators. In
particular, the Riesz transforms on the manifold may not have
Calder\'on-Zygmund kernels, either because one does not have regularity
estimates, or worse because one does not even have size estimates. It
turns out also that the range of $p$ for which one obtains $L^p$
boundedness may not be $(1,\infty)$. See
\cite{ACDH} for a detailed account on all this and Section~\ref{section:applic}
below.

Here, we wish to develop a weighted theory: we want to obtain
weighted $L^p$ estimates for a range of $p$ and for Muckenhoupt
weights with respect to the volume form. Of course, they must
encompass the unweighted estimates so we shall restrict ourselves to
situations where the unweighted theory has been developed. Nothing
new will be done on the unweighted case (except the commutator result
in Section~\ref{section:comm}). We assume that the volume form
is doubling. In that case, we are able to apply
a machinery developed in the first article of our series \cite{AM1}.
In a sense, the results we obtain could have been included in the
latter as an application of the theory there, but we preferred a
separate article to focus here on the geometric aspect of manifolds
and because it also needed technical estimates or ideas from
\cite{ACDH}. As a matter of facts, recent developments in \cite{BZ} allowed us to improve and simplify our result. The reader should have both \cite{ACDH} and
\cite{AM1} handy from now on.

Let $M$ be a complete non-compact  Riemannian manifold with $d$ its
geodesic distance. Let  $\Delta$ be the positive Laplace-Beltrami
operator on $M$ given by
$$
\langle\Delta f, g\rangle = \int_M \nabla f \cdot \nabla g \, d\mu
$$
where $\nabla$ is the Riemannian gradient on $M$ and $\cdot$ is an inner product on $TM$. The Riesz transform is the tangent space valued operator   $\nabla \Delta^{-1/2}$ and it is bounded from $L^2(M,\mu)$ into $L^2(M;TM, \mu)$ by construction.

The manifold $M$ verifies the
doubling volume property if the volume form is doubling:
$$
\mu(B(x,2\,r))\le C\,\mu(B(x,r))<\infty,
$$
for all $x\in M$ and $r>0$ where $B(x,r)=\{y\in M: d(x,y)<r\}$. A Riemannian manifold
$M$ equipped with the geodesic distance and a doubling volume form is a
space of homogeneous type. Non-compactness of $M$  implies  infinite
diameter, which together with the doubling volume property yields
$\mu(M)=\infty$ (see for instance
\cite{Mar}).

One says that the heat kernel $p_t(x,y)$ of the semigroup $e^{-t\Delta}$ has Gaussian upper bounds if  for some constants $c,C>0$ and all $t>0, x,y \in M$,
$$
p_t(x,y) \le  \frac  C{\mu(B(x,\sqrt t))} \, e^{-c \frac{d^2(x,y)}{t}}.
$$
It is known that under doubling  it is a consequence of the same inequality only at  $y=x$  \cite[Theorem 1.1]{G}. We recall the theorem proved in \cite{CD}.

\begin{theor}\label{Rp>2}
Let $M$ be a complete non-compact Riemannian manifold satisfying the
doubling volume property and the Gaussian upper bounds.  Then
$$
\big\|\,|\nabla\Delta^{-1/2}f|\,\big\|_p \le C_p
\|f\|_p \eqno (R_{p})
$$
holds for $1<p< 2$ and all $f$ bounded with compact support.
\end{theor}

 Here,  $|\cdot|$ is the norm on $TM$ associated with the inner product.

We shall set
$$
q_{+}=\sup\big\{p\in (1,\infty)\, : \, (R_p) \ {\rm holds}\big\}.
$$
which  satisfies $q_{+}\ge 2$ under the assumptions of Theorem \ref{Rp>2}. It can
be equal to 2 (\cite{CD}). It is  bigger than 2 assuming further  the stronger $L^2$-Poincar\'e inequalities  (\cite{AC}). It can be  equal to $+\infty$ (see below).

Let us turn to weighted estimates. Properties of Muckenhoupt weights $A_{p}$  and reverse
H\"older classes $RH_{s}$  are reviewed in \cite[Section 2]{AM1}. If
$w\in A_{\infty}(\mu)$, one can define $r_{w}=\inf
\{p>1\, : \, w\in A_{p}(\mu)\}\in [1,\infty)$ and $s_{w}=\sup\{ s>1\,
: \, w\in RH_{s}(\mu)\}\in (1,\infty]$. Given $1\le p_0<q_0\le
\infty$, we introduce the (possibly empty) set
$$
\W_{w}(p_{0},q_{0})
=
\Big(\,p_0\,r_w, \frac{q_0}{(s_w)'}\Big)
=
\big\{
p: p_{0}<p<q_0, w\in A_{\frac{p}{p_0}}(\mu)\cap
RH_{\left(\frac{q_0}{p}\right)'}(\mu)
\big\}.
$$
Here, $q'=\frac {q}{q-1}$ is the conjugate exponent to $q$. And note that $RH_{1}$ means no condition on the weight (besides $A_{\infty}$).

\begin{theor}\label{theor:Rp:w}
Let $M$ be a complete non-compact Riemannian manifold satisfying the
doubling volume property and Gaussian upper bounds. Let $w\in A_{\infty}(\mu)$.
\begin{list}{$(\theenumi)$}{\usecounter{enumi}\leftmargin=.85cm
\labelwidth=.7cm\labelsep=0.25cm \itemsep=0.25cm \topsep=.3cm
\renewcommand{\theenumi}{\roman{enumi}}}

\item For $p\in \W_{w}(1, q_{+})$, the Riesz transform is of strong-type $(p,p)$ with respect to $w\,d\mu$, that is,
\begin{equation}\label{eq:Rp:w}
\big\|\, |\nabla \Delta^{-1/2}f|\, \big\|_{L^p(M,w)}
\le
C_{p,w}\, \|f\|_{L^p(M,w)}
\end{equation}
for all $f$ bounded with compact support.

\item If  $w\in A_{1}(\mu) \cap RH_{(q_{+})'}(\mu) $, then the
Riesz transform is of weak-type $(1,1)$ with respect to $w\,d\mu$,
that is,
\begin{equation}\label{eq:R1:w}
\big\|\, |\nabla \Delta^{-1/2}f|\, \big\|_{L^{1,\infty}(M,w)}
\le
C_{1,w}\, \|f\|_{L^1(M,w)}
\end{equation}
for all $f$ bounded with compact support.
\end{list}

\end{theor}

If   $q_{+}=\infty$  then the Riesz transform is bounded on
$L^p(M,w)$ for $r_{w}<p<\infty$, that is, for $w\in A_{p}(\mu)$, and
we obtain the same weighted theory as for the Riesz transform on
$\RR^n$:

\begin{corol}\label{corol:q+infty}
Let $M$ be a complete non-compact Riemannian manifold satisfying the
doubling volume property and Gaussian upper bounds. Assume that the Riesz
transform has strong type $(p,p)$ with respect to $d\mu$  for all $1<p<\infty$. Then the Riesz
transform has strong type $(p,p)$ with respect to $w\, d\mu$ for all $w\in A_{p}(\mu)$ and
$1<p<\infty$ and it is of weak-type $(1,1)$ with respect to $w\,
d\mu$ for all $w\in A_{1}(\mu)$.
\end{corol}

In \cite[Lemma 4.6]{AM1}, examples of weights in $A_p(\mu)\cap RH_q(\mu)$ are
given. The computations are done in the Euclidean setting, but most
of them can be carried out in spaces of homogeneous type. In
particular, given $f,g\in L^1(M,\mu)$ (or Dirac masses) $1\le r<\infty$ and $1<s\le \infty$, we have that
$w(x)=M_\mu f(x)^{-(r-1)}+ M_\mu g(x)^{1/s}\in A_p(\mu)\cap RH_{q}(\mu)$ ($M_{\mu}$ is the Hardy-Littlewood maximal function) for
all $p>r$ and $q<s$ (and $p=r$ if $r=1$ and $q=s$ if $s=\infty$).
Thus, $r_w\le r$ and $s_w\ge s$.

 We next provide some applications, then proof of our main result  and
eventually we add a short discussion on how to obtain (new) estimates for commutators with bounded mean oscillation functions.

\section{Applications}\label{section:applic}

Unweighted $L^p$ bounds for Riesz transforms in different specific
situations were reobtained in a unified manner in \cite{ACDH} and the
methods used there are precisely those which allowed us to start the
weighted theory. Therefore, it is natural to apply this theory in
return to those situations. Let us concentrate on
five situations (more is done in
\cite{ACDH}). Recall the notion of Poincar\'e inequalities: Let $1\le p<\infty$. One says that $M$ satisfies the $L^p$-Poincar\'e
property, we write $M$ satisfies $(P_{p})$, if there exists $C>0$
such that, for every ball $B$ and every $f$ with $f, \nabla f\in
L^p_{\rm loc}(\mu)$,
 $$
\int_{B}|f-f_{B}|^p\,d\mu\le Cr(B)^p\int_{B} |\nabla f|^p\,d\mu.
\eqno (P_{p})
$$
Here,  $r(B)$ is the radius of $B$, $f_{B}$ is the mean value of $f$
over $B$. It has been proved in \cite{SC1} that  the doubling volume property and $(P_2)$ is equivalent to lower and upper Gaussian estimates on the heat kernel. Recall that $(P_{p})$
implies $(P_{q})$ when $q>p$ (see for instance \cite{HK}).

\subsection{Example without $(P_2)$}

Consider two copies of $\RR^n$ minus the unit ball glued smoothly along their unit circles with $n\ge 2$. It is shown in \cite{CD} that this manifold has  doubling volume form  and Gaussian upper bounds. $(P_2)$ does not hold: in fact, it satisfies $(P_p)$ if and only if $p>n$ (see \cite{HK} in the case of a double-sided cone in $\RR^n$, which is the same). If $n=2$, $(R_p)$  holds if and only if  $p\le 2$ (\cite{CD}). If $n>2$, $(R_p)$ holds if and only if $p<n$ (\cite{CCH}). In any case, we have $q_+=n$. Hence,
$\W_{w}(1,q_{+})= (r_{w}, n /(s_{w})')$ is (contained in) the
range of  $L^p$ boundedness for a given weight provided this is not
empty. In other words, if $1<p<n$ and $w\in A_{p}(\mu) \cap
RH_{(n/p)'}(\mu)$ then one has strong type $(p,p)$ with respect
to $w\,d\mu$. For $p=1$, one has weak-type $(1,1)$ with respect to
$w\, d\mu$ if  $w\in A_{1}(\mu) \cap RH_{n'}(\mu)$.

\subsection{Manifolds with non-negative Ricci curvature}

In this case, the Riesz transform is bounded on (unweighted) $L^p$
for $1<p<\infty$ (\cite{Ba}, \cite{Basurv}). Thus $q_{+}=\infty$. Such
manifolds are known to have  doubling volume form  (see \cite[Theorem 3.10]{Ch}),  $(P_{2})$ and
even $(P_{1})$ \cite{Bu} (see, for instance, \cite{HK} or \cite{parma} for other references). By Corollary \ref{corol:q+infty}, we obtain
strong-type $(p,p)$ for  $1<p<\infty$ and $A_{p}(\mu)$ weights and weak-type (1,1) for
$A_{1}(\mu)$ weights.

\subsection{Co-compact covering manifolds with
polynomial growth deck transformation group} In this case, one has
the doubling volume property and $(P_{2})$ (see
\cite{parma})\footnote{$(P_{1})$ also holds by a discretization
method \cite[Th\'eor\`eme 7.2] {CS} and Poincar\'e inequalities for
discrete groups (see \cite[p.76]{HK}).}. That the Riesz transform is
of  unweighted strong type $(p,p)$ for $1< p\le 2$ is due to
\cite{CD}. For $2<p<\infty$ this is first done in
\cite{ND} and hence $q_+=\infty$.
By Corollary \ref{corol:q+infty}, we obtain
strong-type $(p,p)$ for $1<p<\infty$ and $A_{p}(\mu)$ weights and weak-type (1,1) for
$A_{1}(\mu)$ weights.

\subsection{Conical manifolds with compact basis without boundary}
As mentioned in \cite{ACDH}, this is not strictly speaking a smooth
manifold but it is stochastically complete and this is what is needed
to develop the unweighted theory for the Riesz transform: it is shown in \cite{LHQ} that   $q_{+}$ is a finite value related to the bottom of the spectrum
on the Laplace operator on the compact basis. Also, one has doubling and $(P_{2})$ by \cite{LHQ} and
\cite{CouLi}  (and even $(P_1)$ by using the methods in \cite{GS}).
 Hence,
$\W_{w}(1,q_{+})= (r_{w}, q_{+} /(s_{w})')$ is (contained in) the
range of  $L^p$ boundedness for a given weight provided this is not
empty. In other words, if $1<p<q_{+}$ and $w\in A_{p}(\mu) \cap
RH_{(q_{+}/p)'}(\mu)$ then one has strong type $(p,p)$ with respect
to $w\, d\mu$. For $p=1$, one has weak-type $(1,1)$ with respect to
$w\, d\mu$ if  $w\in A_{1}(\mu) \cap RH_{(q_{+})'}(\mu)$.

\subsection{Lie groups with polynomial volume growth endowed with a sublaplacian}
One starts with left-invariant vector fields $X_{j}$ satisfying the
H\"orman\-der condition and $\mu$ is the left (and right) invariant
Haar measure. The sublaplacian is   $\Delta= - \sum_{j=1}^n X_{j}^2$.
One has the doubling volume property and $(P_{2})$ (and even $(P_1)$)
(see
\cite{V}  or \cite[p. 70]{HK} for a statement and references). The statement of Theorem \ref{theor:Rp:w} applies with no
change to the Riesz transforms $X_{j}\Delta^{-1/2}$. In this case,
$q_{+}=\infty$ from \cite{A}. By Corollary \ref{corol:q+infty}, the
weighted theory for these Riesz transforms is the same as the ones in
$\RR^n$ for $1\le p<\infty$: strong type $(p,p)$ with respect to $w\, d\mu$ holds for $w\in A_p(\mu)$ and
$1<p<\infty$, and weak-type $(1,1)$ with respect to $w\,d\mu$ for
$w\in A_1(\mu)$.

\section{Proof of the main result}

We assume that $M$ satisfies the  doubling volume property  and Gaussian upper bounds.

We first introduce some notation. Given a ball $B$ we set $C_{j}(B)= 4B$ for $j=1$ and $C_{j}(B)= 2^{j+1}\, B \setminus 2^j\, B$ for $j\ge 2$, where
$\lambda\, B$ is the ball co-centered with $B$ and radius $\lambda r(B)$.  We use the notation
$$
\aver{B} h d\mu = \frac {1}{\mu(B)} \int_{B} h\, d\mu,
\qquad\quad
\aver{C_j(B)} h\,d\mu = \frac{1}{\mu(2^{j+1}\, B)} \int_{C_j(B)}
h\,d\mu.
$$
We write $D_\mu$ for the doubling order of $\mu$: $\mu(\lambda B)\le C\,\lambda^{D_\mu}\,\mu(B)$ for every $\lambda>1$. In the sequel $M_\mu$ is the Hardy-Littlewood maximal function with respect to the measure $\mu$ on $M$.

We state a particular case of \cite[Theorem 3.1]{AM1}\footnote{It is stated in the Euclidean setting, see  \cite[Section 5]{AM1} for the extension to spaces of homogeneous type}  that is used in the proof of Theorem \ref{theor:Rp:w}, $(i)$.

\begin{theor}\label{theor:good-lambda:w}
Fix $1<q<\infty$, $a\ge 1$ and $v\in RH_{s'}(\mu)$, $1<s<q$.
Then, there exist $C$ and $K_0\ge 1$ with the
following property: Assume that $F$,  $G$ and $H_1$ are
non-negative measurable functions  on $M$ such that for any ball $B$ there
exist non-negative functions $G_B$ and $H_B$ with $F(x)\le G_B(x)+
H_B(x)$ for a.e. $x\in B$ and, for all $x,\bar{x}\in B$,
\begin{equation}\label{H-Q:G-Q}
\Big(\aver{B} H_B^q\, d\mu \Big)^{\frac1q}
\le
a\, M_\mu F(x)+ H_1(\bar{x}),
\qquad\qquad
\aver{B} G_B\, d\mu
\le
G(x).
\end{equation}
If $1<r\le q/s$ and $F\in L^1(M,\mu)$ (this assumption being only qualitative) we have
\begin{equation}\label{good-lambda:Lp:w}
\|M_\mu F\|_{L^r(M, v)}
\le
C\,\|G\|_{L^r(M, v)}+C\,\|H_1\|_{L^r(M, v)}.
\end{equation}
\end{theor}

\begin{proof}[Proof of  Theorem \ref{theor:Rp:w}, $(i)$]
The argument borrows some ideas from \cite{BZ} which are adapted to the present situation. Fix $w\in A_\infty(\mu)$  and $p\in \W_w(1,q_+)$. Then (see, for instance, \cite[Proposition 2.1]{AM1}) there exist $p_0, q_0$ such that
\begin{equation*}\label{exponents}
1<p_0<p<q_0<q_+
\qquad {\rm and} \qquad
w\in A_{\frac{p}{p_0}}(\mu) \cap RH_{\left(\frac{q_0}{p}\right)'}(\mu).
\end{equation*}
By \cite[Lemma 4.4]{AM1} we have that $v=w^{1-p'}\in A_{p'/q_0'}(\mu)\cap RH_{(p_0'/p')'}(\mu)$. We write $\T=\nabla \Delta^{-1/2}$ and observe that the boundedness of $\T$ from  $L^p(M,w)$ to $L^p(M;TM, w)$ is equivalent to that of $\T^*$  from  $L^{p'}(M; TM, v)$ to $L^{p'}(M,v)$   ---we notice that $\T$ takes scalar valued functions on $M$ to functions valued in  the tangent space (sections) and $\T^*$ the opposite,  $TM$ being equipped with the inner product arising in the definition of the Laplace-Beltrami operator---.

Fix $f\in L^\infty_c (M; TM, \mu)$\footnote{Here and subsequently, the subscript $c$ means with compact support.}, and  write $h=\T^* f$ and $F=|h|^{q_0'}$. Notice that $F\in L^1(M, \mu)$ since $\T^*$ is bounded from $L^{q_0'}(M;TM,\mu)$ to $L^{q_0'}(M,\mu)$ as $1<q_0<q_+$ and $\T$ is bounded from $L^{q_0}(M,\mu)$ to $L^{q_0}(M;TM,\mu)$. We pick $\A_{r}=I- (I-e^{-r^2\, \Delta})^m$ with $m$ large enough. Given a ball $B$ we write $r_B$ for its radius. Then,
$$
F
\le
G_B+H_B
\equiv
2^{q_0'-1}\,|(I-\A_{r_B})^* h|^{q_0'}
+
2^{q_0'-1}\,|\A_{r_B}^* h|^{q_0'}.
$$

We first estimate $H_B$. Set $q=p_0'/q_0'$ and observe that by duality there exists $g\in L^{p_0}(B,d\mu/\mu(B))$ with norm $1$ such that for all $x\in B$
\begin{align*}
\Big(\aver{B} H_B^q\,d\mu\Big)^{\frac1{q\,q_0'}}
&
\lesssim
\mu(B)^{-1}\,\int_{M} |h|\,|\A_{r_B} g|\,d\mu
\\
&\lesssim
\sum_{j=1}^\infty 2^{j\,D_\mu} \Big(\aver{C_j(B)} |h|^{q_0'}\,d\mu\Big)^\frac1{q_0'}\,\Big(\aver{C_j(B)} |\A_{r_B} g|^{q_0}\,d\mu\Big)^\frac1{q_0}
\\
&\le
M_\mu F(x)^\frac1{q_0'}\sum_{j=1}^\infty 2^{j\,D_\mu} \Big(\aver{C_j(B)} |\A_{r_B} g|^{q_0}\,d\mu\Big)^\frac1{q_0}.
\end{align*}
To estimate the summands we use the Gaussian upper bound  on $p_{t}(x,y)$, so that for any fixed integer $m$   there exist $c,C>0$
such that   for all $j\ge 1$,   all ball $B$, all $g \in
L^1(M,\mu)$ supported in $B$ and all $1\le k
\le m$,
\begin{equation}\label{eq:L1Linfty}
\sup_{{C_j(B)}} |e^{- k\, r_B^2\, \Delta}g|
\le
C\, e^{-c\,4^j}\, \aver{B} |g|\,d\mu.
\end{equation}
Then expanding $\A_{r}=I-
(I-e^{-r^2\, \Delta})^m$  we conclude that
\begin{equation}\label{est-HB}
\Big(\aver{B} H_B^q\,d\mu\Big)^{\frac1{q\,q_0'}}
\lesssim
M_\mu F(x)^\frac1{q_0'}\sum_{j=1}^\infty 2^{j\,D_\mu}\,e^{-c\,4^j}\,
\Big(\aver{B} |g|^{q_0}\,d\mu\Big)^\frac1{q_0}
\lesssim
M_\mu F(x)^\frac1{q_0'}.
\end{equation}

We next estimate $G_B$. Using duality there exists $g\in L^{q_0}(B,d\mu/\mu(B))$ with norm $1$ such that for all $x\in B$
\begin{align*}
\Big(\aver{B} G_B\,d\mu\Big)^{\frac1{q_0'}}
&
\lesssim
\mu(B)^{-1}\,\int_{M} |f|\,|\T(I-\A_{r_B}) g|\,d\mu
\\
&\lesssim
\sum_{j=1}^\infty 2^{j\,D_\mu} \Big(\aver{C_j(B)} |f|^{q_0'}\,d\mu\Big)^\frac1{q_0'}\,\Big(\aver{C_j(B)} |\T(I-\A_{r_B}) g|^{q_0}\,d\mu\Big)^\frac1{q_0}
\\
&\le
M_\mu(|f|^{q_0'})(x)^\frac1{q_0'}\sum_{j=1}^\infty 2^{j\,D_\mu} \Big(\aver{C_j(B)} |\T(I-\A_{r_B}) g|^{q_0}\,d\mu\Big)^\frac1{q_0}.
\end{align*}
 To estimate the terms in the sum we use the following auxiliary result whose proof is given below.

\begin{lemma}\label{lemma:beta}
For all $\beta \in [1, \tilde q_{+}) \cup [1,2]$, one has the
following estimate:
 for all $m\ge 1$, there exists $C>0$ such that for
all $j\ge 2$, all ball $B$, all $g \in L^1(M,\mu)$ with support
in $B$,
\begin{equation}\label{eq:beta}
\Big(
\aver{C_j(B)} |\nabla \Delta^{-1/2} (I-e^{-r(B)^2\,
\Delta})^m g|^{\beta}\,d\mu\Big)^{\frac1{\beta}}
\le
C4^{-jm}\,\aver{B} |g|\,d\mu.
\end{equation}
\end{lemma}

Here,  $\tilde q_{+}$ is defined as
the supremum of those $p\in (1,\infty)$ such that for all $t>0$,
  \begin{equation}
\label{eq:Gp} \big\|\, |\nabla  e^{-t\, \Delta}f  |\, \big\|_{p} \le
C\,t^{-1/2} \|f\|_{p}.
\end{equation}
By analyticity of the heat semigroup, one always have $\tilde
q_{+}\ge q_{+}$. Under the doubling volume property and $(P_{2})$, it is
shown in
\cite[Theorem 1.3]{ACDH} that $q_{+}=\tilde q_{+}$. We do not know if the equality holds or not under doubling and Gaussian upper bounds.

This lemma allows us to conclude right away with the terms where $j\ge 2$. When  $j=1$, we use that $\T$ is bounded from $L^{q_0}(M,\mu)$ to $L^{q_0}(M; TM, \mu)$ as $1<q_0<q_+$. Also, applying \eqref{eq:L1Linfty} it follows easily that
\begin{align}
\aver{4\,B} |\T(I-\A_{r_B}) g|^{q_0}\,d\mu
&\lesssim
\frac1{\mu(4\,B)}\,\Big(\int_{B}|g|^{q_0}\,d\mu+\sum_{j=1}^\infty \int_{C_j(B)} |\A_{r_B} g|^{q_0}\,d\mu\Big)
\nonumber
\\
&\lesssim
\aver{B} |g|^{q_0}\,d\mu\ \sum_{j=1}^\infty 2^{j\,D_{\mu}}\,e^{-c\,4^j}
\lesssim
\aver{B} |g|^{q_0}\,d\mu.
\label{est-GB:j=1}
\end{align}
Using this and \eqref{eq:beta} (with $\beta=q_0< q_+ \le \tilde{q}_+$) we conclude the estimate for $G_B$:
\begin{align}
\Big(\aver{B} G_B\,d\mu\Big)^{\frac1{q_0'}}
&
\lesssim
M_\mu(|f|^{q_0'})(x)^\frac1{q_0'}\sum_{j=1}^\infty 2^{j\,D_\mu} 4^{-j\,m} \Big(\aver{B} |g|^{q_0}\,d\mu\Big)^\frac1{q_0}
\nonumber
\\
&\le
C\,M_\mu(|f|^{q_0'})(x)^\frac1{q_0'}
=
G(x)^{\frac1{q_0'}},
\label{est-GB}
\end{align}
provided $m>D_\mu/2$.
With these estimates in hand we can use Theorem \ref{theor:good-lambda:w} with $r=p'/q_0'$, $q= p_0'/q_0'$ and $H_1\equiv 0$. Notice that  $v\in RH_{s'}(\mu)$ with $s=p_0'/p'$,  $1<s<q<\infty$ and $r=q/s$. Hence, using $v \in A_r(\mu)$ we obtain the desired estimate
\begin{equation}\label{est-T*}
\|\T^* f\|_{L^{p'}(M,v)}^{q_0'}
\le
\|M_\mu F\|_{L^r(M,v)}
\lesssim
\|M_\mu(|f|^{q_0'})\|_{L^r(M,v)}
\lesssim
\||f|\|_{L^{p'}(M,v)}^{q_0'}.
\end{equation}
\end{proof}

\begin{proof} [Proof of  Lemma \ref{lemma:beta}]
First, this estimate is known for $\beta=2$ (see \cite{ACDH}). Also,
the inequality for a fixed $\beta_0$ implies the same one for all
$\beta$ with $1\le\beta \le \beta_0$. It suffices to treat the case
$\beta>2$, which happens only if $\tilde q_{+}>2$.

We use  a trick from \cite[Proof of Lemma 3.1]{ACDH}. Fix a ball
$B$, with radius $r$, and $f\in L^\infty(M, \mu)$ supported in $B$. We
have
$$
\nabla \Delta^{-1/2} (I-e^{-r^2\Delta})^m f = \int_0^\infty
g_{r}(t) \, \nabla e^{-t\, \Delta} f\,  {dt}
$$
where $g_{r}\colon \RR^+ \to \RR$ is a function such that
\begin{equation}
\label{eq:gr} \int_0^\infty |g_{r}(t)|\,e^{-\frac {c\,  4^j\, r^2}
t}\, \frac {dt} {\sqrt t}  \le C_{m}\,  4^{-j\, m}.
\end{equation}
By definition of $\tilde q_+$ and the argument of
\cite[p. 944]{ACDH} we have
$$
\Big(\int_M | \nabla_x\, p_t(x,y)|^\beta \,
e^{\gamma \frac{ d^2(x,y)}{t}} \, d\mu(x)\Big)^{1/\beta}
\le \frac C{ {\sqrt t}
\left[\mu(B(y,\sqrt{t}))\right]^{1-1/\beta}},
$$
for all $t>0$ and $y\in M$,   with  $\gamma>0$ depending on $\beta$.
 This implies that for all $j\ge 2$, $y\in B$ and
all $t>0$,
$$
\Big(
\aver{C_j(B)} |\nabla_{x}\,
p_{t}(x,y)|^{\beta}\,d\mu(x)\Big)^{1/{\beta}}
\lesssim \frac 1{\sqrt t}\, e^{-\frac
{c\,  4^j\, r^2} t} \frac 1 {\mu(B(y,\sqrt t))^{1-1/\beta}\,
\mu(2^{j+1}\, B)^{1/\beta}}.
$$
Using the doubling property, $\mu(2^{j+1}\, B) \sim \mu (B(y,
2^{j+1}\, r))$ uniformly in $y\in B$ and
$$
\frac {\mu(B(y, 2^{j+1}\, r)) }{\mu(B(y,\sqrt t)) }
\lesssim \max \Big\{ 1, \frac {2^j\, r}{\sqrt t}\, \Big\}^{D_\mu}.
$$
Hence, with another
$c>0$,
$$
\Big(
\aver{C_j(B)} |\nabla_{x}\,
p_{t}(x,y)|^{\beta}\,d\mu(x)\Big)^{1/{\beta}}
\lesssim
\frac 1{\sqrt t}\, e^{-\frac {c\,  4^j\, r^2} t} \frac 1 {
\mu(2^{j+1}\, B)}
\le
\frac 1{\sqrt t}\, e^{-\frac {c\,  4^j\, r^2} t} \frac 1 { \mu(B)}
.
$$
We conclude using Minkowski's integral inequality and \eqref{eq:gr}
that the left hand side of \eqref{eq:beta} is bounded by
\begin{align*}
\lefteqn{\hskip-2cm
\int_{0}^\infty |g_{r}(t)| \int_B |f(y)|\, \Big( \aver{C_j(B)}
|\nabla_{x}\, p_{t}(x,y)|^{\beta}\,d\mu(x)\Big)^{1/{\beta}}
d\mu(y)\, dt}
\\
&  \lesssim
\int_{0}^\infty |g_{r}(t)|  \frac 1{\sqrt t}\, e^{-\frac {c\,
4^j\, r^2} t} \, dt \,   \aver{B} |f|\, d\mu
\lesssim
4^{-j\,m}\, \aver{B} |f|\,d\mu .
\end{align*}
\end{proof}

The following result, used to prove Theorem \ref{theor:Rp:w}, $(ii)$,  is taken from \cite[Theorem 8.8 \& Remark 8.10]{AM1} (see also \cite[Section 8.4]{AM1} for the extension to spaces of homogeneous type)

\begin{theor}\label{theor:B-K:small:w}
Let $1\le p_0<q_0\le \infty$ and  $w\in A_{\infty}(\mu)$. Let $T$ be a sublinear operator defined
on  $L^2(M,\mu)$ and $\{\A_r\}_{r>0}$ be a family of operators acting
from $L_c^\infty(M,\mu)$  into $L^2(M,\mu)$. Assume the following conditions:
\begin{list}{$(\theenumi)$}{\usecounter{enumi}\leftmargin=1cm
\labelwidth=0.7cm\itemsep=0.1cm\topsep=.3cm
\renewcommand{\theenumi}{\alph{enumi}}}

\item There exists  $q\in \W_{w}(p_{0},q_{0})$ such that  $T$ is bounded from $L^q(M,w)$ to $L^{q,\infty}(M,w)$.

\item For all $j\ge 1$, there exist a constant $\alpha_{j}$ such that   for any ball $B$ with $r(B)$ its radius and
for any $f\in L_c^\infty(M, \mu)$ supported in $B$,
\begin{equation}\label{small:A:w}
\Big(
\aver{C_j(B)} |\A_{r(B)}f|^{q_0}\,d\mu\Big)^{\frac1{q_0}}
\le
\alpha_j\,\Big(\aver{B} |f|^{p_0}\,d\mu\Big)^{\frac1{p_0}}.
\end{equation}

\item There exists $\beta>(s_{w})'$, \textit{i.e.} $w\in RH_{\beta'}(\mu)$,
with the following property: for all $j\ge 2$, there exist a constant
$\alpha_{j}$ such that  for any ball $B$ with $r(B)$ its radius and
for any $f\in L_c^\infty(M,\mu)$ supported in $B$ and for  $j\ge 2$,
\begin{equation}\label{small:T:I-A:w}
\Big( \aver{C_j(B)} |T(I-\A_{r(B)})f|^{\beta}\,d\mu\Big)^{\frac1\beta}
\le
\alpha_j\,\Big(\aver{B} |f|^{p_0}\,d\mu\Big)^{\frac1{p_0}}.
\end{equation}

\item $\sum_j \alpha_j\,2^{D_{w}\,j}<\infty$ for $\alpha_{j}$ in $(b)$ and $(c)$, where $D_w$ is the doubling constant of $w\, d\mu$.
\end{list}
If $w\in A_1(\mu)\cap RH_{(q_0/p_0)'}(\mu)$ then $T$ is of weak-type $(p_0,p_0)$ with respect to $w\, d\mu$, that is, for all $f\in L^\infty_c(M, \mu)$,
$$
\|T f\|_{L^{p_0,\infty}(M,w)} \le C\,  \| f\|_{L^{p_0}(M,w)}.
$$
\end{theor}

\begin{proof}[Proof of Theorem \ref{theor:Rp:w},  $(ii)$]
We are going to apply Theorem \ref{theor:B-K:small:w} with $p_0=1$ and $q_0=q_+$. Thus we need to check that the four items hold.
Fix $w\in A_{1}(\mu) \cap RH_{(q_{+})'}(\mu) $.  By  \cite[Proposition 2.1]{AM1}, there exists $1<q<q_+$ such that $w\in A_q(\mu)\cap RH_{(q_+/q)'}(\mu)$. This means that $q\in \W_w(1,q_+)$, therefore by $(i)$, $T=|\nabla\Delta^{-1/2}|$ is bounded on $L^q(M,w)$ and so $(a)$ holds.

We pick $\A_{r}=I- (I-e^{-r^2\, \Delta})^m$ with $m$ large enough to be chosen.  Notice that expanding $\A_r$, \eqref{eq:L1Linfty} yields $(b)$ with $\alpha_j=C\,e^{-c\,4^j}$. To see $(c)$ we apply Lemma \ref{lemma:beta} with $(s_w)'<\beta$ ---notice that such $\beta$ exist: we have $q_+ \le \tilde q_+$ and $w\in RH_{(q_+)'}(\mu)$ implies $(s_w)'<q_+$---. Then, we obtain $(c)$ with $\alpha_j=C\,4^{-j\,m}$. Finally, we pick $m>D_w/2$ so that $(d)$ holds and therefore Theorem \ref{theor:B-K:small:w} gives the weak-type $(1,1)$ with respect to $w\, d\mu$.
\end{proof}

\section{Commutators}\label{section:comm}

Let us write again $\T=\nabla\Delta^{-1/2}$ and take $b\in \BMO(M,\mu)$ (the
space of bounded mean oscillation functions on $M$). We define the first
order commutator  $\T_b^1 g=[b,\T]g=b\,\T g - \T(b\,g)$, and for $k\ge 2$
the $k$-th order commutator is $\T_b^k =[b,\T_b^{k-1}]$. Here $g,b$ are scalar valued and
$\T_b^k g$ is valued in the tangent space.

\begin{theor} Under the assumptions of Theorem \ref{theor:Rp:w},
 $\T_b^k$ satisfies \eqref{eq:Rp:w} for each  $k\ge 1$, that is, it is bounded from $L^p(M,w)$  into $L^p(M;TM,w)$ under the same conditions on $w,p$.
\end{theor}

This theorem applies in particular to the five situations described in Section
\ref{section:applic}. Note that even the unweighted $L^p$ estimates
for the commutators are new.

\begin{proof}
The proof is similar to that of Theorem \ref{theor:Rp:w} using again the ideas in \cite{BZ} and we point out the main changes. We only consider the case $k=1$: the general case follows by induction and  the details are left to the reader (see \cite[Section 6.2]{AM1} for similar arguments). As in \cite[Lemma 6.1]{AM1} it suffices to assume qualitatively  $b\in L^\infty(M,\mu)$ and  quantitatively $\|b\|_{\BMO(M,\mu)} = 1$ and get uniform bounds.

We proceed as before working with $\T_b^1$ in place of $\T$. Write $F=|(\T_b^1)^* f|^{q_0'}$ with $f\in L^\infty_c(M;TM,\mu)$ and observe that
$F\in L^1(M,\mu)$ as $b\in L^\infty(M,\mu)$ and $\T^*$ is bounded from $L^{q_0'}(M;TM,\mu)$ into $L^{q_0'}(M,\mu)$ (as $1<q_0<q_+$) ---we observe that this is the only place where we use that $b\in L^\infty(M,\mu)$---. Fixing $B$ we write $\hat{b}=b-b_B$ and
decompose $\T_b^1$ as $\T_b^1 g=
-\T (\hat{b} g)
+
\hat{b}\,\T g
$.
Using this equality one sees $(\T_b^1)^*=-(\T^*)_b^1$. Then we have
\begin{align*}
F
&=
|(\T_b^1)^* f|^{q_0'}
=
|(\T^*)_b^1 f|^{q_0'}
\le
2^{q_0'-1}\,|\hat{b}\,\T^*f|^{q_0'}
+
2^{q_0'-1}\,|\T^*(\hat{b}\,f)|^{q_0'}
\\
&\le
\big(
2^{q_0'-1}\,|\hat{b}\,\T^*f|^{q_0'}
+
4^{q_0'-1}\,|(I-\A_{r_B})^*\T^*(\hat{b}\,f)|^{q_0'}\big)
+
4^{q_0'-1}\,|\A_{r_B}^*\T^*(\hat{b}\,f)|^{q_0'}
\\
&=
G_B+H_B.
\end{align*}
We estimate $H_B$. By duality we pick $g$ as before and obtain
\begin{align*}
&\Big(\aver{B} H_B^q\,d\mu\Big)^{\frac1{q\,q_0'}}
=
C\,
\mu(B)^{-1}\,\int_{M} \T^* (\hat{b}\,f)\,\A_{r_B} g\,d\mu
\\
&\qquad\qquad
=
C\,
\mu(B)^{-1}\,\int_{M} \big(-(\T^*)_b^1 f + \hat{b}\,\T^* f\big)\,\A_{r_B} g\,d\mu
\\
&
\qquad\qquad
\lesssim
\mu(B)^{-1}\,\int_{M} |(\T^*)_b^1 f|\,|\A_{r_B} g|\,d\mu
+
\mu(B)^{-1}\,\int_{M} |\hat{b}|\,|\T^* f|\,|\A_{r_B} g|\,d\mu
=
I+II.
\end{align*}
The estimate for $I$ follows as in \eqref{est-HB} by using \eqref{eq:L1Linfty}: for all $x\in B$,
$$
I
\lesssim
\sum_{j=1}^\infty 2^{j\,D_\mu} \Big(\aver{C_j(B)} F\,d\mu\Big)^\frac1{q_0'}\,\Big(\aver{C_j(B)} |\A_{r_B} g|^{q_0}\,d\mu\Big)^\frac1{q_0}
\lesssim
M_\mu F(x)^\frac1{q_0'}.
$$
We pick $1<\tilde{q}<\infty$ and use H\"older's inequality to obtain that for all $\bar{x}\in B$
\begin{align*}
II
&
\lesssim
\sum_{j=1}^\infty 2^{j\,D_\mu}
\Big(\aver{C_j(B)} |\T^*f|^{q_0'}\,d\mu\Big)^\frac1{q_0'}\,
\Big(\aver{C_j(B)} |\A_{r_B} g|^{q_0\,s}\,d\mu\Big)^\frac1{q_0\,s}\,
\Big(\aver{2^{j+1}\,B} |\hat{b}|^{q_0\,s'}\,d\mu\Big)^\frac1{q_0\,s'}
\\
&\lesssim
\|b\|_{\BMO(M,\mu)}
M_\mu(|\T^*f|^{q_0'})(\bar{x})^\frac1{q_0'}
\sum_{j=1}^\infty 2^{j\,D_\mu} e^{-c\,4^j}\,(1+j)
\lesssim
M_\mu(|\T^*f|^{q_0'})(\bar{x})^\frac1{q_0'},
\end{align*}
where we have used \eqref{eq:L1Linfty} and John-Nirenberg's inequality. Collecting $I$ and $II$ we conclude the first estimate in \eqref{H-Q:G-Q} with  $H_1=M_\mu(|\T^*f|^{q_0'})$.
Let us  write $G_{B,1}$ and $G_{B,2}$ for each of the terms that define  $G_B$ and we estimate them in turn. Take $\delta>1$ to be chosen and use John-Nirenberg's inequality: for any $x\in B$ we have
\begin{align*}
\aver{B} G_{B,1}\,d\mu
&=
C\,\aver{B}
|\hat{b}\,\T^*f|^{q_0'}\,d\mu
\lesssim
\Big(\aver{B} |\T^*f|^{q_0'\,\delta}\,d\mu\Big)^\frac1{\delta}\,
\Big(\aver{B} |b-b_B|^{q_0'\,\delta'}\,d\mu\Big)^\frac1{\delta'}
\\
&\lesssim \|b\|_{\BMO(M,\mu)}^{q_0'}
M_\mu(|\T^*f|^{q_0'\,\delta})(\bar{x})^\frac1{\delta}.
\end{align*}
To estimate $G_{B,2}$ we proceed as with $G_B$ in the proof of Theorem \ref{theor:Rp:w}. Let $g$ be the corresponding dual function and use again John-Nirenberg's inequality
\begin{align*}
\Big(\aver{B} G_{B,2}\,d\mu\Big)^\frac1{q_0'}
&=
C\,
\Big(\aver{B}
|(I-\A_{r_B})^*\,\T^*(\hat{b}f)|^{q_0'}\,d\mu\Big)^\frac1{q_0'}
\\
&\lesssim
\mu(B)^{-1}\,\int_{M} |\hat{b}|\,|f|\,|\T(I-\A_{r_B}) g|\,d\mu
\\
&\lesssim
\sum_{j=1}^\infty 2^{j\,D_\mu}
\Big(\aver{2^{j+1}\,B} |b-b_B|^{q_0'\,\delta'}\,d\mu\Big)^\frac1{q_0'\delta'}
\,
\Big(\aver{C_j(B)}  |f|^{q_0'\,\delta}\,d\mu\Big)^\frac1{q_0'\,\delta}\,
\\
&\hskip4cm
\Big(\aver{C_j(B)} |\T(I-\A_{r_B}) g|^{q_0}\,d\mu\Big)^\frac1{q_0}
\\
&\lesssim
M_\mu(|f|^{\delta\,q_0'})(x)^\frac1{\delta\,q_0'}\sum_{j=1}^\infty 2^{j\,D_\mu}\,(1+j) \Big(\aver{C_j(B)} |\T(I-\A_{r_B}) g|^{q_0}\,d\mu\Big)^\frac1{q_0}
\\
&\lesssim
M_\mu(|f|^{\delta\,q_0'})(x)^\frac1{\delta\,q_0'}.
\end{align*}
where in the last estimate we have proceeded as in \eqref{est-GB} using \eqref{eq:beta} for $j\ge 2$ and \eqref{est-GB:j=1} for $j=1$. Gathering what has been obtained for $G_{B,1}$ and $G_{B,2}$ we conclude the second  estimate in \eqref{H-Q:G-Q} with  $G=
M_\mu(|\T^*f|^{q_0'\,\delta})^\frac1{\delta}
+
M_\mu(|f|^{\delta\,q_0'})^\frac1{\delta}
$.

We apply Theorem \ref{theor:good-lambda:w} as in the proof of Theorem \ref{theor:Rp:w}. In this case, we observe that as $v\in A_r(\mu)$, there exists $1<\delta<r$ so that $v\in A_{r/\delta}(\mu)$. Then,
the desired estimate follows
\begin{align*}
\|(\T_b^1)^* f\|_{L^{p'}(M,v)}^{q_0'}
&
\le
\|M_\mu F\|_{L^r(M,v)}
\lesssim
\|G\|_{L^r(M,v)}+ \|H_1\|_{L^r(M,v)}
\\
&\le
\|M_\mu(|f|^{\delta\,q_0'})^\frac1{\delta}\|_{L^r(M,v)}
+
\|M_\mu(|\T^*f|^{\delta\,q_0'})^\frac1{\delta}\|_{L^r(M,v)}
\\
&\lesssim
\|f\|_{L^{p'}(M;TM,v)}^{q_0'}
+
\|\T^*f\|_{L^{p'}(M,v)}^{q_0'}
\lesssim
\|f\|_{L^{p'}(M;TM,v)}^{q_0'}
\end{align*}
where we have used that $\T^*$ is bounded from $L^{p'}(M;TM,v)$ into $L^{p'}(M,v)$ (see \eqref{est-T*}).
\end{proof}

\end{document}